\documentclass[12pt]{article}
\usepackage[cp1251]{inputenc}
\usepackage[english]{babel}
\usepackage{amssymb,amsmath,amsthm}
\usepackage{graphicx, epsfig, color, hyperref}
\newtheorem{predl}{Proposition}
\usepackage{caption}

\begin{document}
\title{On exact representations of the motions group of Galilean plane}
\author{Dmitry Efimov, Igor Kostyakov, Vasiliy Kuratov}
\maketitle
\begin{abstract}
The Pimenov algebra with two generators is defined and some of its properties
are shown. Some exact matrix over the Pimenov algebra representations of the
motions group of Galilean plane (the Galilean group) are considered. A
geometric interpretation of them is giving. We consider also a exact
representation of the Galilean group by elements of Grassmann algebra.
\end{abstract}

\section{Introduction}
In this paper we describe several exact matrix and one hypercomplex representation of the group motions of the Galilean plane,
which arises in many mathematical as well as theoretical physics problems (see, e.g., \cite{Vil}, \cite{Jag1}).

The Galilean plane is one of the nine possible two-dimensional spaces with constant curvature \cite{Jag1}.
A distance $d$ between points $M_1(x_1,y_1)$ and $M_2(x_2,y_2)$ here is
\begin{eqnarray}\label{f_1}
d=\begin{cases}
|x_1-x_2|,&\text{$x_1\not=x_2$;}\\
|y_1-y_2|,&\text{$x_1=x_2$}.
\end{cases}
\end{eqnarray}
By {\it motion} of Galilean plane is understood its transformation
preserves the distance between the points and does not change orientation.
It is known that any such movement can be represented as the composition translations
along the coordinate axes and rotation around the origin (the Galilean boost).
By definition of distance (\ref{f_1}), rotations represent irregular plane translations  along the vertical axis.
The set of all motions of the Galilean plane is a group with respect to the motions composition.
Denote it by $G(2)$.

From the point of view of physics Galilean plane is the simplest model of the nonrelativistic
spacetime. 
One of the coordinates is interpreted as a space coordinate,
and the other as time, and plane motions are Galilean transformations
together with translations of reference system.
Because of its relatively simple structure, it is a good training model
to test and refine different theoretical constructs.
Group $G(2)$ is considered also often in applications to quantum mechanics.
In this case, it is called usually {\it Heisenberg group}.
What is said above  explains the continued interest in this group in theoretical physics literature,
including recent years.

The paper is organized as follows.
In the second section we define and give some properties of real Pimenov algebra with two generators and its subalgebra --- the algebra of dual numbers, on the use which a further material is based.
In the next two sections  exact $3\times 3$ and $2\times 2$ matrix representations of the group $G(2)$ are considered.
In the last section the representation of $G(2)$ with elements of Grassmann algebra is considered
At the end of paper a application is given, which collect, for convenience, all considered
in this paper exact representations of the group motions of Galilean plane.

\section{The Pimenov algebra and dual numbers}

By {\it the Pimenov algebra} with two generators $D_2({\mathbf R})$ (or simply $D_2$) we will mean the associative algebra with the unit
and two generators $\iota_1$, $\iota_2$ (dual units) such that
\begin{equation}\label{fp1}
\iota_1^2=\iota_2^2=0,\ \ \iota_1\iota_2=\iota_2\iota_1\not=0.
\end{equation}
Therefore any element of  the algebra  $D_2$ is represented uniquely in the form:
\begin{equation}\label{fp2}
a=a_0+a_1\iota_1+a_2\iota_2+a_3\iota_1\iota_2,\ \ a_i\in{\mathbf R}. 
\end{equation}
By analogy with complex numbers $a_0$ and $a-a_0$ are called {\it the real} and {\it the imaginary part} of $a$ respectively.
Addition, multiplication, multiplication by scalar operations are defined in naturally manner in the algebra $D_2$.
They are commutative and associative. And the division operation is determined partially in contrast to complex numbers.
Namely the inverse element of the element (\ref{fp2}) has the form:
$$
a^{-1}=\frac{1}{a_0^2}\left[a_0-a_1\iota_1-a_2\iota_2+\left(2\frac{a_1a_2}{a_0}-a_3\right)\iota_1\iota_2\right].
$$
From this formula it follows that the division operation is defined only on elements with nonzero real part $(a_0\not=0)$.
For arbitrary element (\ref{fp2}) one can introduce the conjugation by the generator $\iota_2$:
$$
\tilde{a}=a_0+a_1\iota_1-a_2\iota_2-a_3\iota_1\iota_2.
$$
It is not hard to prove that  $\widetilde{ac}=\tilde{a}\tilde{c}$, $a,c\in D_2$. 
The elements of the Pimenov algebra can be assigned to a class of so called {\it hypercomplex numbers} \cite{KS}.

Any matrix $A$ over the algebra $D_2$ can be represented uniquely in the form:
$$
A=A_0+\iota_1A_1+\iota_2A_2+\iota_1\iota_2A_3, 
$$
where $A_i$ are matrices of the same size as $A$ with elements from ${\mathbf R}$.
Matrices $A_0$ and $A-A_0$ are called  {\it the real} and {\it the imaginary part} of the matrix $A$ respectively.
Let $M_n(D_2)$ be the set of all $n\times n$ matrices with elements from $D_2$.
Denote also the set of all $n\times n$ matrices with real elements by $M_n({\mathbf R})$.
Than not hard to show that the following proposition holds.
\begin{predl}\label{p0}
A matrix  $A=A_0+\iota_1 A_1+\iota_2A_2+\iota_1\iota_2A_3$ from $M_n(D_2)$ is invertible if and only if
its real part $A_0$ is invertible in $M_n({\mathbf R})$. In this case the inverse matrix has the form:
$$
A^{-1}=A_0^{-1}\left[A_0-\iota_1A_1-\iota_2A_2+\iota_1\iota_2
 \left(A_1A_0^{-1}A_2+A_2A_0^{-1}A_1-A_3\right)\right]A_0^{-1}.
$$
\end{predl}
The set of all invertible matrices from $M_n(D_2)$ forms a group under multiplication. Denote it by $GL_n(D_2)$.
Call a matrix from $M_n(D_2)$ {\it nondegenerate} if its determinant is invertible in $D_2$.
In view of previous proposition and the fact that the real part of the product of elements from $D_2$ 
is equal to the product of real parts, we obtain the following propositon.
\begin{predl}\label{p0_1}
The real part of the determinant of a matrix $A\in M_n(D_2)$ is equal to the determinant of the real part of  $A$:
$$
\mathrm{Re}\det A=\det\mathrm{Re}\, A. 
$$
A matrix $A\in M_n(D_2)$ is invertible if and only if it is nondegenerate.
\end{predl}

Consider any matrix $A\in M_n(D_2)$. Denote by $\tilde{A}$ the matrix,
which  is obtained from $A$ by application of the conjugation by $\iota_2$ to all its elements. 
Introduce for matrices from $M_n(D_2)$ the  {\it dual conjugation} $\star$ by the rule:
$$
A^{\star}\stackrel{def}=\tilde{A}^T.
$$
The set of matrices from $GL_n(D_2)$ such that $A^\star=A^{-1}$ and $\det{A}=1$ forms a group, 
which denote by $SU(D_2)$ by analogy with the special unitary group.

Functions of $D_2$-arguments are defined by their Taylor expansion
in which the imaginary part of an argument plays a role of the increment \cite{Dim}. 
Thus taking into account properties of dual units, we obtain
$$
f(a_0+a_1\iota_1+a_2\iota_2+a_3\iota_1\iota_2)=f(a_0)+f'(a_0)(a_1\iota_1+a_2\iota_2+a_3\iota_1\iota_2)+f''(a_0)a_1a_2\iota_1\iota_2.
$$
For example,
$$
e^{a_1\iota_1+a_2\iota_2+a_3\iota_1\iota_2}=1+a_1\iota_1+a_2\iota_2+(a_3+a_1a_2)\iota_1\iota_2.
$$

If we put $a_2=a_3=0$ in the expansion (\ref{fp2}), then we get elements of the form
$$
a=a_0+\iota a_1,\ \ \iota^2=0,
$$
which are called  {\it the dual numbers} \cite{Dim}, \cite{Jag}.
The set of all dual numbers forms a subalgebra in the Pimenov algebra $D_2$
and hence all the above obviously carries over to them.

\section{Exact $3\times 3$ matrix representations of the group $G(2)$}
One of the most popular representation of this group has the form:
\begin{eqnarray}\label{f1}
g=\left(
\begin{array}{ccc}
1&0&0\\
a&1&0\\
b&\theta&1
\end{array}\right),\ \ a,b,\theta\in{\mathbf R}.
\end{eqnarray}	
Here the parameters $a$ and $b$ characterize translations along the coordinate axes, the parameter $\theta$ --- a rotation around the origin.
If to each point $(x,y)$ of  Galilean plane assign the column vector $(1,x,y)^T$, then
the representation (\ref{f1}) acts on Galilean plane by the rule:
\begin{eqnarray}\label{f2}
\left(
\begin{array}{ccc}
1&0&0\\
a&1&0\\
b&\theta&1
\end{array}\right)
\left(
\begin{array}{c}
1\\
x\\
y
\end{array}\right)=
\left(
\begin{array}{c}
1\\
x+a\\
y+\theta x+b
\end{array}\right).
\end{eqnarray}
In kinematical interpretation the variable $x$ is the time, the variable $y$ is the space coordinate,
the parameter $\theta$ defines the value of the constant velocity of the moving inertial system
relatively fixed.
Each element (\ref{f1}) can be uniquely represented as the  product of elements of the one-parametric subgroups:
$$
g=\left(
\begin{array}{ccc}
1&0&0\\
a&1&0\\
0&0&1
\end{array}\right)
\left(
\begin{array}{ccc}
1&0&0\\
0&1&0\\
b&0&1
\end{array}\right)
\left(
\begin{array}{ccc}
1&0&0\\
0&1&0\\
0&\theta&1
\end{array}\right).
$$
The infinitesimal operators of this one-parametric subgroups in the unit 
$$
A_1=\left(
\begin{array}{ccc}
0&0&0\\
1&0&0\\
0&0&0
\end{array}\right),\ \ \ 
A_2=\left(
\begin{array}{ccc}
0&0&0\\
0&0&0\\
1&0&0
\end{array}\right),\ \ \ 
A_3=\left(
\begin{array}{ccc}
0&0&0\\
0&0&0\\
0&1&0
\end{array}\right).
$$
with commuting relations 
\begin{eqnarray}\label{f2_3}
[A_1,A_2]=0,\ \ \ [A_2,A_3]=0,\ \ \ [A_3,A_1]=A_2
\end{eqnarray}
generate the Lie algebra $g(2)$ of the motion group of Galilean plane.

Galilean plane belongs to a family of nine possible two-dimensional spaces with constant curvature.
It is known that all these spaces in a unified manner are modeled on the connected components of the unit sphere in
three-dimensional spaces  ${\mathbf R}^3(j_1,j_2)$, consisting of vectors $v=(x,j_1y,j_1j_2z)$, where $x,y,z\in{\mathbf R}$,
and $j_1,j_2$ take the values $1,i,\iota_k$; $\iota_k^2=0$, $\iota_1\iota_2=\iota_2\iota_1\not=0$, $k=1,2$ \cite{Pim}.
Accordingly, the group of motions of these nine two-dimensional spaces with constant curvature
locally isomorphic to the rotations groups $SO(3,j_1,j_2)$ of the spaces ${\mathbf R}^3(j_1,j_2)$ \cite{G-97}.
In a certain basis  $SO(3,j_1,j_2)$ groups consist of matrices of the form
\begin{eqnarray}\label{f3}
A=\left(
\begin{array}{ccc}
a_{11}&j_1a_{12}&j_1j_2a_{13}\\
j_1a_{21}&a_{22}&j_2a_{23}\\
j_1j_2a_{31}&j_2a_{32}&a_{33}
\end{array}\right),\ A^TA=AA^T=E,\ \det{A}=1,
\end{eqnarray}
where $a_{ij}\in{\mathbf R}$.
Galilean plane correspond to the values of parameters $j_1=\iota_1$, $j_2=\iota_2$.
It is easy to see that the general element of $SO(3,\iota_1,\iota_2)$ in a certain parameterization has the form
\begin{eqnarray}\label{f4}
g=\left(
\begin{array}{ccc}
\sigma_1&-\iota_1\sigma_1\sigma_2 a&-\iota_1\iota_2 (\sigma_2b-a\theta)\\
\iota_1 a&\sigma_2&-\iota_2\sigma_1\theta\\
\iota_1\iota_2 b&\iota_2\theta&\sigma_1\sigma_2
\end{array}\right),
\end{eqnarray}
where $\sigma_i=\pm 1$, $a,b,\theta\in {\mathbf R}$.
The group $SO(3,\iota_1,\iota_2)$  is not simply connected, it consists of four connected components, which correspond to four combinations of parameters 
$\sigma_1,\sigma_2$.
Locally, as mentioned above, the group $SO(3,\iota_1,\iota_2)$ is isomorphic to the motions group of Galilean plane $G(2)$.
Namely, if in the formula (\ref{f4}) put $\sigma_1=\sigma_2=1$, we obtain the general element of $G(2)$:
\begin{eqnarray}\label{f5}
g=\left(
\begin{array}{ccc}
1&-\iota_1 a&-\iota_1\iota_2 (b-a\theta)\\
\iota_1 a&1&-\iota_2\theta\\
\iota_1\iota_2 b&\iota_2\theta&1
\end{array}\right).
\end{eqnarray}
Indeed, every element (\ref{f5}) can be uniquely represented as the product
\begin{eqnarray}\label{f5_1}
g=\left(
\begin{array}{ccc}
1&-\iota_1 a&0\\
\iota_1 a&1&0\\
0&0&1
\end{array}\right)
\left(
\begin{array}{ccc}
1&0&-\iota_1\iota_2 b\\
0&1&0\\
\iota_1\iota_2 b&0&1
\end{array}\right)
\left(
\begin{array}{ccc}
1&0&0\\
0&1&-\iota_2\theta\\
0&\iota_2\theta&1
\end{array}\right)
\end{eqnarray}
of the one-parameter subgroups elements corresponding to the rotations on the finite angle in the coordinate
planes of the space ${\mathbf R}^3(\iota_1,\iota_2)$.
Infinitesimal operators of these one-parameter subgroups in the unit are:
$$
A_1=\left(
\begin{array}{ccc}
0&-\iota_1&0\\
\iota_1&0&0\\
0&0&0
\end{array}\right),\ 
A_2=\left(
\begin{array}{ccc}
0&0&-\iota_1\iota_2\\
0&0&0\\
\iota_1\iota_2&0&0
\end{array}\right),\ 
A_3=\left(
\begin{array}{ccc}
0&0&0\\
0&0&-\iota_2\\
0&\iota_2&0
\end{array}\right).
$$
Their commutation relations obviously coincide with the commutation relations (\ref{f2_3}).

Galilean plane itself corresponds to the connected component of the unit sphere
$x^2+(\iota_1 y)^2+(\iota_1\iota_2 z)^2=1$ or  $x^2=1$ of the space ${\mathbf R}^3(\iota_1,\iota_2)$,
i.e. in this case just to the plane $x=1$.
The transformation (\ref{f5}) act on it by the rule:
\begin{eqnarray}\label{f6}
\left(
\begin{array}{ccc}
1&-\iota_1 a&-\iota_1\iota_2 (b-a\theta)\\
\iota_1 a&1&-\iota_2\theta\\
\iota_1\iota_2 b&\iota_2\theta&1
\end{array}\right)
\left(
\begin{array}{c}
1\\
\iota_1 y\\
\iota_1\iota_2 z
\end{array}\right)
=\left(
\begin{array}{c}
1\\
\iota_1(y+a)\\
\iota_1\iota_2  (z+\theta y+b)
\end{array}\right),
\end{eqnarray}
(cf. (\ref{f2})). 
Thus the formula (\ref{f5}) gives another matrix representation of $G(2)$
with orthogonal matrices. But in this case we need to go beyond
the field of real numbers and consider matrices  over Pimenov algebra $D_2({\mathbf R})$.

\section{Exact $2\times 2$ matrix representations of the group $G(2)$}
The group $G(2)$ has also exact $2\times 2$ matrix representations.
Consider some of them. The following proposition holds.
\begin{predl}\label{p1}
The group $G(2)$ is isomorphic to the three-parameter subgroup $G$ of matrices 
\begin{eqnarray}\label{f7}
g=\left(
\begin{array}{cc}
e^{\iota_2\phi}&\iota_1(\beta+\iota_2\gamma)\\
-\iota_1(\beta-\iota_2\gamma) &e^{-\iota_2\phi}
\end{array}\right),\ \phi,\beta,\gamma\in \mathbf{R},\ \iota_k^2=0,\ \iota_1\iota_2=\iota_2\iota_1
\end{eqnarray}
of the group $SU(D_2)$. 
\end{predl}

\begin{proof}
One-parameter subgroups of $G$ have the form:
$$
g(a)=\left(
\begin{array}{cc}
1&\iota_1\frac{a}{2}\\
-\iota_1\frac{a}{2}&1
\end{array}\right),\ 
g(b)=\left(
\begin{array}{cc}
1&\iota_1\iota_2\frac{b}{2}\\
\iota_1\iota_2\frac{b}{2}&1
\end{array}\right),\
g(\theta)=\left(
\begin{array}{cc}
e^{\tfrac{\iota_2 \theta}{2}}&0\\
0&e^{-\tfrac{\iota_2 \theta}{2}}
\end{array}\right).
$$
Their infinitesimal operators in the unit are:
$$
A_1=\frac{1}{2}\left(
\begin{array}{cc}
0&\iota_1\\
-\iota_1&0
\end{array}\right),\ \ \  
A_2=\frac{1}{2}\left(
\begin{array}{cc}
0&\iota_1\iota_2\\
\iota_1\iota_2&0
\end{array}\right),\ \ \ 
A_3=\frac{1}{2}\left(
\begin{array}{cc}
\iota_2&0\\
0&-\iota_2
\end{array}\right).
$$
Obviously, their commutation relations coincide with the commutation relations (\ref{f2_3})
of the Lie algebra of $G(2)$.
It is also easy to check that $g^\star=g^{-1}$ and $\det g=1$ (see Section 2), i.e. $g\in SU(D_2)$.
\end{proof} 
We now consider the action of this representation on Galilean plane.
It is closely linked with the action of the previous representation.
Namely, we associate to each matrix (\ref{f5_1})  the matrix
\begin{eqnarray}\label{f10}
\begin{aligned}
&u=\left(
\begin{array}{cl}
e^{\tfrac{\iota_2\theta}{2}}&\iota_1\left[\frac{a}{2}+\iota_2\frac{b}{2}\right] e^{-\tfrac{\iota_2\theta}{2}}\\
-\iota_1\left[\frac{a}{2}-\iota_2\frac{b}{2}\right] e^{\tfrac{\iota_2\theta}{2}}&\ \ \ \ \ e^{-\tfrac{\iota_2\theta}{2}}
\end{array}\right)=\\
&=\left(
\begin{array}{cc}
1&\iota_1\frac{a}{2}\\
-\iota_1\frac{a}{2}&1
\end{array}\right)
\left(
\begin{array}{cc}
1&\iota_1\iota_2\frac{b}{2}\\
\iota_1\iota_2\frac{b}{2}&1
\end{array}\right)
\left(
\begin{array}{cc}
e^{\tfrac{\iota_2\theta}{2}}&0\\
0&e^{-\tfrac{\iota_2\theta}{2}}
\end{array}\right).
\end{aligned}
\end{eqnarray}
of the form (\ref{f7}).  It is easy to verify that this correspondence is a group isomorphism.
By definition, the dual-conjugate matrix to the matrix $u$ has the form:
$$
u^\star=\left(
\begin{array}{cl}
e^{-\tfrac{\iota_2\theta}{2}}&-\iota_1\left[\frac{a}{2}+\iota_2\frac{b}{2}\right] e^{-\tfrac{\iota_2\theta}{2}}\\
\iota_1\left[\frac{a}{2}-\iota_2\frac{b}{2}\right] e^{\tfrac{\iota_2\theta}{2}}&\ \ \ \ \ \ \ \ \  e^{\tfrac{\iota_2\theta}{2}}
\end{array}\right).
$$
Put each vector-column $v=(1,\iota_1 y,\iota_1\iota_2 z)^T$ into correspondence with  the matrix
\begin{eqnarray}\label{f11}
h_v=\left(
\begin{array}{cc}
0&\iota_1 \left[\frac{y}{2}+\iota_2\frac{z}{2}\right])\\
\iota_1\left[\frac{y}{2}-\iota_2\frac{z}{2}\right])&1
\end{array}\right).
\end{eqnarray}
Then the following action of this representation in the space of matrices (\ref{f11}) corresponds to rotation (\ref{f6}):
\begin{eqnarray}\label{f12}
\begin{aligned}
uh_vu^\star=&\left(
\begin{array}{cc}
0&\iota_1\left[\frac{y'}{2}+\iota_2\frac{z'}{2}\right]\\
\iota_1\left[\frac{y'}{2}-\iota_2\frac{z'}{2}\right]&1
\end{array}\right)=h_{v'},\\
&y'=y+a,\ \ z'=z+\theta y+b
\end{aligned}
\end{eqnarray}.

\begin{figure}[!ht]
\begin{center}
\includegraphics[scale=0.5]{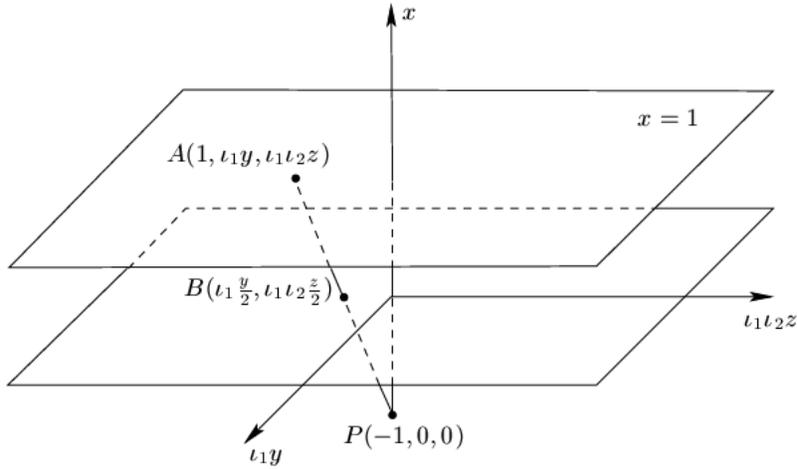}
\caption{ The stereographic projection in the space ${\mathbf R}^3(\iota_1,\iota_2)$}\label{pic1}
\end{center}
\end{figure}

To give a geometric interpretation of this correspondence we consider in the space ${\mathbf R}^3(\iota_1,\iota_2)$
the stereographic projection with  a pole in $P(-1,0,0)$.
It put each point $A(1,\iota_1 y,\iota_1\iota_2 z)$ of the connected component of the unit sphere $x^2=1$ 
into correspondence with the point $B(\iota_1 \frac{y}{2},\iota_1\iota_2\frac{z}{2})$ of the plane $x=0$ (Fig.~\ref{pic1}).
In turn  one can associate the point $B$ uniquely with hypercomplex numbers
\begin{equation}\label{f13}
\xi=\iota_1\left(\frac{y}{2}+\iota_2\frac{z}{2}\right),\ \ \ \ \ \widetilde{\xi}=\iota_1\left(\frac{y}{2}-\iota_2\frac{z}{2}\right).
\end{equation}
Consider the rotation (\ref{f12}) of the connected component of the unit sphere.
Suppose that the point  $A(1,\iota_1 y,\iota_1\iota_2 z)$ of the sphere
transforms  by the rotation to the point $A'(1,\iota_1 y',\iota_1\iota_2 z')$ and
the stereographic projection takes $A'$ to the point $B'(\iota_1 \frac{y'}{2},\iota_1\iota_2\frac{z'}{2})$.
Suppose also that $B'$ associate by the rule (\ref{f13}) with hypercomplex numbers
$$
\eta=\iota_1\left(\frac{y'}{2}+\iota_2\frac{z'}{2}\right),\ \ \ \ \ \widetilde{\eta}=\iota_1\left(\frac{y'}{2}-\iota_2\frac{z'}{2}\right).
$$
It is easy to show that hypercomplex numbers $\xi$, $\eta$, $\widetilde{\xi}$ and $\widetilde{\eta}$ connected to each other the following relations:
\begin{eqnarray}\label{f15}
\eta=
\frac{e^{\tfrac{\iota_2\theta}{2}}\xi+\iota_1\left(\frac{a}{2}+\iota_2\frac{b}{2}\right) e^{-\tfrac{\iota_2\theta}{2}}}
{-\iota_1\left(\frac{a}{2}-\iota_2\frac{b}{2}\right) e^{\tfrac{\iota_2\theta}{2}}\xi+e^{-\tfrac{\iota_2\theta}{2}}},\ \ 
\widetilde{\eta}=
\frac{e^{-\tfrac{\iota_2\theta}{2}}\widetilde{\xi}+\iota_1\left(\frac{a}{2}-\iota_2\frac{b}{2}\right) e^{\tfrac{\iota_2\theta}{2}}}
{-\iota_1\left(\frac{a}{2}+\iota_2\frac{b}{2}\right) e^{-\tfrac{\iota_2\theta}{2}}\widetilde{\xi}+e^{\tfrac{\iota_2\theta}{2}}}. 
\end{eqnarray} 
If we introduce on the plane $x=0$ of the space ${\mathbf R}^3(\iota_1,\iota_2)$ homogeneous coordinates
$(\iota_1\left[y_1+\iota_2 z_1\right],y_2+\iota_2 z_2)^T$, $y_1, y_2, z_1, z_2\in {\mathbf R}$, $y_2\not=0$
such that $\xi\thicksim (\iota_1\xi_1,\xi_2)^T$, $\eta\thicksim (\iota_1\eta_1,\eta_2)^T$, i.e. $\xi=\frac{\iota_1\xi_1}{\xi_2}$,  $\eta=\frac{\iota_1\eta_1}{\eta_2}$, then equalities (\ref{f15}) one can rewrite as:
$$
\left(
\begin{array}{c}
\iota_1\eta_1\\
\eta_2
\end{array}
\right)=
\left(
\begin{array}{cl}
e^{\tfrac{\iota_2\theta}{2}}&\iota_1\left[\frac{a}{2}+\iota_2\frac{b}{2}\right] e^{-\tfrac{\iota_2\theta}{2}}\\
-\iota_1\left[\frac{a}{2}-\iota_2\frac{b}{2}\right] e^{\tfrac{\iota_2\theta}{2}}&\ \ \ \ \ e^{-\tfrac{\iota_2\theta}{2}}
\end{array}\right)
\left(
\begin{array}{c}
\iota_1\xi_1\\
\xi_2
\end{array}
\right),\ 
$$
$$
\left(
\begin{array}{c}
\iota_1\widetilde{\eta_1}\\
\widetilde{\eta_2}
\end{array}
\right)^T=
\left(
\begin{array}{c}
\iota_1\widetilde{\xi_1}\\
\widetilde{\xi_2}
\end{array}
\right)^T
\left(
\begin{array}{cl}
e^{-\tfrac{\iota_2\theta}{2}}&-\iota_1\left[\frac{a}{2}+\iota_2\frac{b}{2}\right] e^{-\tfrac{\iota_2\theta}{2}}\\
\iota_1\left[\frac{a}{2}-\iota_2\frac{b}{2}\right] e^{\tfrac{\iota_2\theta}{2}}&\ \ \ \ \ \ \ \ \  e^{\tfrac{\iota_2\theta}{2}}
\end{array}\right),\ 
$$
or briefly
\begin{eqnarray}\label{f16}
\eta=u\xi,\ \ \ \ \ \eta^\star=\xi^\star u^\star.
\end{eqnarray}
Now normalize the homogeneous coordinates so that $y_2 + \iota_2z_2 = 1$. Then, it is easy to see that
$$
\xi\xi^\star=
\left(
\begin{array}{cc}
0&\iota_1 \left[\frac{y}{2}+\iota_2\frac{z}{2}\right]\\
\iota_1\left[\frac{y}{2}-\iota_2\frac{z}{2}\right]&1
\end{array}\right)=h_v,\ \ \ \eta\eta^\star=h_{v'}.
$$
Hence, multiplying the equality (\ref{f16}), we obtain the equality (\ref{f12}).
Thus, the rotation (\ref{f12}) of the connected component of the unit sphere in the space ${\mathbf R}^3(\iota_1,\iota_2)$ corresponds
to a fractional-linear transformation of the special form (\ref{f15}) of the plane $x=0$, in which the unit sphere appears by the stereographic projection,
and this correspondence is a group homomorphism. 

In the paper \cite{Mac} are discussed also in detail two-dimensional spaces of constant curvature and their groups of motions. 
Adopted here is similar to the approach in \cite{Pim}, but some of it is different.
As a result, among other things, one more exact $2\times 2$ matrix representation of the motions group of Galilean plane is obtained there.
And it is shown that the group $G(2)$ is isomorphic to the three-parameter group $G$ of matrices
\begin{eqnarray}\label{f19_1}
g=\left(
\begin{array}{cc}
e^{\iota\phi}&\zeta+\iota\eta\\
0 &e^{-\iota\phi}
\end{array}\right),\ \ \ \phi,\zeta,\eta\in \mathbf{R},\ \iota^2=0.
\end{eqnarray}
If we assign the matrix 
$$
w=\left(
\begin{array}{cc}
e^{\tfrac{\iota\theta}{2}}&\frac{a+\iota b}{2}e^{-\tfrac{\iota\theta}{2}} \\
0&e^{-\tfrac{\iota\theta}{2}}
\end{array}\right)
$$
of the form (\ref{f19_1}) to each matrix (\ref{f5})
and assign the matrix
\begin{eqnarray}\label{f19_3}
p_v=\left(
\begin{array}{cc}
-1&y+\iota z\\
0&1
\end{array}\right),
\end{eqnarray}
to each vector-column $v=(1,\iota_1 y,\iota_1\iota_2 z)^T$, then the action
$$
wp_vw^{-1}=\left(
\begin{array}{cc}
-1&y'+\iota z'\\
0&1
\end{array}\right)=p_{v'};\ \ 
y'=y+a,\ z'=z+\theta y+b
$$
of this representation in the space of matrices (\ref{f19_3}) correspond to transformation (\ref{f6}).

At the conclusion of this section  we give another convenient $2\times 2$ matrix representation of the group $G(2)$.
Put each element (\ref{f1}) of the group $G(2)$ into correspondence with the three-parameter matrix
\begin{eqnarray}\label{f20}
g=\left(
\begin{array}{cc}
e^{\iota\theta}&a+\iota b\\
0&1
\end{array}\right)=
\left(
\begin{array}{cc}
1+\iota\theta&a+\iota b\\
0&1
\end{array}\right),\ \ \ \iota^2=0.
\end{eqnarray}
It is not hard to show that this  mapping is a group isomorphism.
If one put  each point of the Galilean plane
with coordinates $(x,y)$ into correspondence with column vector $(x+\iota y,1)^T$,
then the transformation (\ref{f20}) will be act on the plane by the rule: 
$$
\left(
\begin{array}{cc}
e^{\iota\theta}&a+\iota b\\
0&1
\end{array}\right)
\left(
\begin{array}{c}
x+\iota y\\
1
\end{array}\right)=
\left(
\begin{array}{c}
x+a+\iota(y+\theta x+b)\\
1
\end{array}\right).
$$
Obviously, this expression is the motion of the Galilean plane coinciding with the motion (\ref{f2}).

\section{The hypercomplex representation of $G(2)$}
At the conclusion give another, the so-called {\itshape hypercomplex} representation of the group $G(2)$.
Note that it can be seen as part of a broader scheme that uses the so-called {\it generalized quaternions} \cite{Mac2}. 
Consider the Grassmann algebra $\Lambda({\mathbf R^2})$, i.e. an associative real algebra with unit and with two generators 
$e_1$, $e_2$ satisfying relations
\begin{equation}\label{f42}
e_1^2=e_2^2=0,\ e_1e_2=-e_2e_1.
\end{equation}
The general element of the algebra $\Lambda({\mathbf R^2})$ can be written as
\begin{equation}\label{f40}
q=\alpha_0+\alpha_1e_1+\alpha_2e_2+\alpha_3e_1e_2,
\end{equation}
where  $\alpha_i\in {\mathbf R}$.
Introduce the conjugation and the norm in $\Lambda({\mathbf R^2})$  by the rule:
$$
\bar{q}=\alpha_0-\alpha_1e_1-\alpha_2e_2-\alpha_3e_1e_2,\ \ \ |q|^2=q\bar{q}=\alpha_0^2.
$$
Denote by $\Lambda^1({\mathbf R^2})$ the set of all elements of the algebra $\Lambda({\mathbf R^2})$  with unit norm and  with $\alpha_0=1$.
Not hard to see that it form a group under multiplication. There is the following interesting fact:
\begin{predl}
The Galilean group $G(2)$ is isomorphic to the group of all elements of the Grassmann algebra $\Lambda({\mathbf R^2})$
with unit scalar part:
$$
G(2)\cong \Lambda^1({\mathbf R^2}).
$$
\end{predl}
\begin{proof}
The general element (\ref{f7}) of the group $G(2)$ can be uniquely
represented as follows:
\begin{equation}\label{f44}
\begin{aligned}
&\left(
\begin{array}{cc}
e^{\iota_2\phi}&\iota_1(\beta+\iota_2\gamma)\\
-\iota_1(\beta-\iota_2\gamma)&e^{-\iota_2\phi}
\end{array}\right)=
\left(
\begin{array}{cc}
1&0\\
0&1
\end{array}
\right)+
\phi\left(
\begin{array}{cc}
\iota_2&0\\
0&-\iota_2
\end{array}
\right)+\\
&+\beta
\left(
\begin{array}{cc}
0&\iota_1\\
-\iota_1&0
\end{array}
\right)+
\gamma\left(
\begin{array}{cc}
0&\iota_1\iota_2\\
\iota_1\iota_2&0
\end{array}
\right)=E+\phi E_1+\beta E_2+\gamma E_1E_2.
\end{aligned}
\end{equation}
It is easy to see that the matrices  $E_1,E_2$ satisfy relations (\ref{f42})
and the coefficient of $E$ is $1$.
This proves the isomorphism $G(2)$ and $\Lambda^1({\mathbf R^2})$.
\end{proof}
Consider one of the actions of the group $G(2)$ in this representation.
For this we take a real Clifford algebra $Cl_3((\mathbf R))$ with three generators $e_1$, $e_2$, $e_3$
satisfying relations
\begin{equation}\label{f44_1}
e_ie_j+e_je_i=0,\ i\not=j;\ \ e_1^2=e_2^2=0,\ \ e_3^2=1.
\end{equation}
One can give this algebra a matrix interpretation, if we consider the matrices $E_1$, $E_2$ and
$$
E_3=\left(
\begin{array}{cc}
-1&0\\
0&1
\end{array}\right)
$$
as generators.
Obviously, the group $\Lambda^1({\mathbf R^2})$ is a subgroup
of the multiplicative group of the algebra $Cl_3({\mathbf R})$.
Note now that one can associate uniquely  to each column vector $v=(1,\iota_1 y,\iota_1\iota_2 z)^T$ 
    the matrix
\begin{equation}\label{f45}
\begin{aligned}
\left(
\begin{array}{cc}
-1&\iota_1\left(y+\iota_2z\right)\\
\iota_1\left(y-\iota_2z\right)&1
\end{array}\right)&
=\left(
\begin{array}{cc}
-1&0\\
0&1
\end{array}
\right)
+y
\left(
\begin{array}{cc}
0&\iota_1\\
\iota_1&0
\end{array}
\right)+\\
+z
\left(
\begin{array}{cc}
0&\iota_1\iota_2\\
-\iota_1\iota_2&0
\end{array}
\right)&=
\left(E+yE_2+zE_1E_2\right)E_3
\end{aligned}
\end{equation}
instead of the matrix (\ref{f11}) in the formula (\ref{f12}).
This expansion shows that the set of matrices (\ref{f45}) is a subspace of the algebra $Cl_3({\mathbf R})$.
Suppose now that an element $q$ of $\Lambda^1({\mathbf R^2})$ corresponds to an element $u$ of the group $G(2)$
by the rule (\ref{f44}).
Then not hard to see that the element $u^\star$ will corresponds to the element $\bar{q}$
and the rotation (\ref{f12}) can be written as follows Clifford form:
$$
 qq_v\bar{q}=q_{v'},
$$
where $q_v$ and $q_{v'}$ are elements of the algebra $Cl_3({\mathbf R})$ which correspond to matrices $h_v$ and $h_{v'}$
by the rule (\ref{f45}).

\vspace{3mm}
The work has been partially supported by RFBR grant 08-01-90010 --- Belarus and the program
``Mathematical problems of nonlinear dynamics'' of the Presidium of Russian Academy of Sciences.

\appendix
\section{Exact matrix and hypercomplex representations of the motions group of Galilean plain}
\noindent
{\bfseries 1)} The ``standart'' representation:

\begin{eqnarray*}
g=\left(
\begin{array}{ccc}
1&0&0\\
a&1&0\\
b&\theta&1
\end{array}\right),\ \ a,b,\theta\in{\mathbf R}.
\end{eqnarray*}

\noindent
{\bfseries 2)} The representation with orthogonal matrices over the algebra $D_2({\mathbf R})$: 

$$
g=\left(
\begin{array}{ccc}
1&-\iota_1 a&-\iota_1\iota_2 (b-a\theta)\\
\iota_1 a&1&-\iota_2\theta\\
\iota_1\iota_2 b&\iota_2\theta&1
\end{array}\right),\ a,b,\theta\in \mathbf{R},\ \iota_k^2=0,\ \iota_1\iota_2=\iota_2\iota_1.
$$

\noindent
{\bfseries 3)} The representation with unimodular matrices over the algebra $D_2({\mathbf R})$:

\begin{eqnarray*}
g=\left(
\begin{array}{cc}
e^{\iota_2\phi}&\iota_1(\beta+\iota_2\gamma)\\
-\iota_1(\beta-\iota_2\gamma) &e^{-\iota_2\phi}
\end{array}\right),\ \ \phi,\beta,\gamma\in \mathbf{R},\ \iota_k^2=0,\ \iota_1\iota_2=\iota_2\iota_1.
\end{eqnarray*}

\noindent
{\bfseries 4)} The $2\times 2$ representation with unimodular matrices with dual elements:
\begin{eqnarray*}
g=\left(
\begin{array}{cc}
e^{\iota\phi}&\zeta+\iota\eta\\
0&e^{-\iota\phi}
\end{array}\right),\ \ \phi,\zeta,\eta\in \mathbf{R},\ \iota^2=0.
\end{eqnarray*}

\noindent
{\bfseries 5)} The $2\times 2$ representation with dual elements:
\begin{eqnarray*}
g=\left(
\begin{array}{cc}
e^{\iota\theta}&a+\iota b \\
0&1
\end{array}\right),\ \ a,b,\theta\in \mathbf{R},\ \iota^2=0.
\end{eqnarray*}

\noindent
{\bfseries 6)} The representation with elements of the Grassmann algebra:
$$
g=1+\alpha_1e_1+\alpha_2e_2+\alpha_3e_1e_2,\ \ \alpha_i\in {\mathbf R},\ \ e_1^2=e_2^2=0,\ \ e_1e_2=-e_2e_1.
$$


\begin{thebibliography}{100}
\bibitem{Vil}
{\bf Vilenkin N.Ya., Klimyk A.U.}~Representation of Lie groups and special functions. Volume 2. Kluwer Academic Publishers, 1992.
\bibitem{Dim}
{\bf Dimentberg F.M.}~The screw calculus and its applications in mechanics. Foreign Technology Division, Wright-Patterson Air Force Base, Ohio, 1968. 155 p.
\bibitem{KS}
{\bf Kantor I.L., Solodovnikov A.S., Shenitzer A.}~Hypercomplex numbers: an elementary introduction to algebras. Springer, 1989. 
\bibitem{Pim}
{\bf Pimenov R.I.}~Unified axiomatics of spaces with the maximum group of motions. {\it Litovski Matematicheski Sbornik}. 1965. V.5. P. 457-486 (in Russian).
\bibitem{Jag}
{\bf Yaglom I.M.}~Complex numbers in geometry. Academic Press, 1968. 243 p.
\bibitem{Jag1}
{\bf Yaglom I.M.}~A simple non-Eucledian geometry and its physical basis: an elementary account of Galilean geometry and the
Galilean principle of relativity. New-York: Springer-Verlag, 1979. (Translated from the Russian).
\bibitem{G-97}
{\bf Gromov N.A., Kostyakov I.V., Kuratov V.V.}~FRT quantization theory for the nonsemisimple Cayley-Klein groups. arXiv:q-alg/9711024
\bibitem{Mac}
{\bf McRae A.}~Clifford algebras and possible kinematics. {\it Symmetry, integrability and geometry: methods and applications.} 2007. V.3. 29 p.
arXiv:0707.2869
\bibitem{Mac2}
{\bf McRae A.}~Clifford fibrations and  possible kinematics. {\it Symmetry, integrability and geometry: methods and applications.} 2009. V.5. 18 p.
arXiv:0907.2394

\end{thebibliography}
\end{document}